\documentclass{article}
\usepackage{amsmath}
\usepackage{amssymb}
\usepackage{amsthm}

\newcommand{\Q}{\mathbb{Q}}
\newcommand{\Qx}{\Q [x]}
\newcommand{\Z}{\mathbb{Z}}

\newcommand{\lra}{\longrightarrow}

\newcommand{\nxn}{(n \times n)}
\newcommand{\w}{w}

\newtheorem{proposition}{Proposition}

\newtheorem{theorem}{Theorem}
\newtheorem{lemma}{Lemma}

\usepackage[all,knot,poly]{xy}

\begin{document}
\title{Proof of a Conjectured Formula for the Alexander Invariant}
\author{Peter Lee}
\maketitle

\begin{abstract}
In this paper we prove the validity of a formula for computing the Alexander invariant, which was originally conjectured by Bar-Natan and Dancso in \cite{DBN-WKO}.
\end{abstract}

\tableofcontents

\section{Introduction}

In recent work \cite{DBN-WKO}, Bar-Natan and Dancso conjecture that a certain formula which they describe computes the Alexander invariant of a knot.  The procedure involves writing down a certain matrix based on data read off from the crossings in the knot, and computing its determinant (which we denote $\beta(K)(x)$).  Extensive computer checking reported in \cite{DBN-WKO} confirmed the validity of Bar-Natan and Dancso's conjecture for prime knots of up to 11 crossings.  The purpose of this paper is to prove that this conjectured formula does indeed give the Alexander invariant.

In his paper \cite{Alex}, Alexander gave a procedure to write down a certain matrix (not obviously related to Bar-Natan's matrix) with entries in a polynomial ring $\Qx$ for every knot, and showed that the determinant $\Delta(K)(x)$ of this matrix is an invariant of the knot (up to a factor of $\pm x^l$ for $l\in \Z$).  The strategy of this paper is to show that the formula described by Bar-Natan and Dancso in \cite{DBN-WKO} implicitly determines a series of row and column operations which convert Alexander's matrix into a matrix whose determinant is $\beta(K)(x)$ (up to signs and powers of $x$).

In Section 2 of this paper we first explain the formula given in \cite{DBN-WKO}, and then recall Alexander's original procedure for computing $\Delta(K)(x)$.  Then in Section 3 we give our proof that these two procedures give the same invariant.

\subsection{Terminological Conventions}
\label{Conventions}
Throughout the paper we will let $K$ be a long knot represented by an oriented diagram with $n$ crossings, numbered from $1$ to $n$ in any order.  The arcs of the long knot $K$ divide the plane into $n+2$ regions.  Two of these regions are unbounded, while the other $n$ regions are bounded by the arcs of the knot.  We further number these bounded regions from $1$ to $n$ in any order.

\subsection{Acknowledgements}

Many thanks to Dror Bar-Natan and Zsuzsi Dancso, to whom the conjecture proven here is due, and more particularly to the former, for reading and providing very useful comments on a draft of this paper.

\section{Background}

\subsection{The Invariant of Bar-Natan and Dancso}

We now describe the procedure given in \cite{DBN-WKO}.

If $K$ is a long knot, with numbered crossings and regions as per our convention (see Section \ref{Conventions}), and one imagines an ant walking along the knot in the direction of the orientation, then for $1\leq i\leq n$, the `span' $Sp(i)$ of crossing $i$ is defined to be the (oriented) open segment of the knot between the two times the ant passes through crossing $i$.  We form an $\nxn$ matrix $T$ by the rule:

\begin{itemize}
\item $T_{ij}=1$ if the ant crosses \emph{under} crossing $j$ at any point as it walks along $Sp(i)$; and
\item $T_{ij}=0$ otherwise.
\end{itemize}

For instance, the following picture shows a long knot with the span of crossing 2 marked (i.e., the dotted line between the two stars *), along with the $T_{2,j}$.  Vertex numbers are shown, but region numbers are omitted.

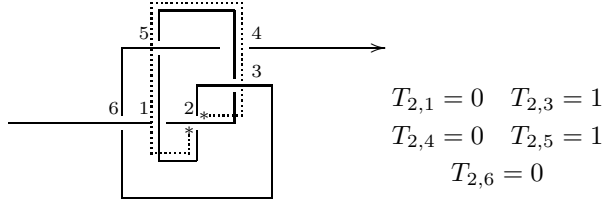
\[
\xy
(-20,0)*{
\xy
(-30,0)*{} = "1";
(-11,0)*{} = "2";
(-9,0)*{} = "3";
(0,0)*{} = "4";
(-10,-5)*{} = "5";
(-10,15)*{} = "6";
(0,15)*{} = "7";
(0,6)*{} = "8";
(0,4)*{} = "9";
(-5,-5)*{} = "11";
(-5,-1)*{} = "12";
(-5,1)*{} = "13";
(-5,5)*{} = "14";
(5,5)*{} = "15";
(-4,1)*{\scriptstyle{*}} = "16";
(1,1)*{} = "17";
(1,16)*{} = "18";
(-11,16)*{} = "19";
(-11,-4)*{} = "20";
(-6,-4)*{} = "21";
(-6,-1)*{\scriptstyle{*}} = "22";
(-12,2)*{\scriptstyle{1}};
(-6,2)*{\scriptstyle{2}};
(3,7)*{\scriptstyle{3}};
(3,12)*{\scriptstyle{4}};
(-12,12)*{\scriptstyle{5}};
(-16,2)*{\scriptstyle{6}};
(5,5)*{} = "23";
(5,-10)*{} = "24";
(-15,-10)*{} = "25";
(-15,-1)*{} = "26";
(-15,1)*{} = "27";
(-15,10)*{} = "28";
(-10,9)*{} = "29";
(-10,11)*{} = "30";
(-2,10)*{} = "31";
(2,10)*{} = "32";
(20,10)*{} = "33";
{\ar@{-} "1"; "2"};
{\ar@{-} "3"; "4"};
{\ar@{-} "4"; "9"};
{\ar@{-} "8"; "7"};
{\ar@{-} "7"; "6"};
{\ar@{-} "6"; "30"};
{\ar@{-} "5"; "11"};
{\ar@{-} "11"; "12"};
{\ar@{-} "13"; "14"};
{\ar@{-} "14"; "15"};
{\ar@{.} "16"; "17"};
{\ar@{.} "17"; "18"};
{\ar@{.} "18"; "19"};
{\ar@{.} "19"; "20"};
{\ar@{.} "20"; "21"};
{\ar@{.} "21"; "22"};
{\ar@{-} "15"; "24"};
{\ar@{-} "24"; "25"};
{\ar@{-} "25"; "26"};
{\ar@{-} "27"; "28"};
{\ar@{-} "28"; "31"};
{\ar@{-} "29"; "5"};
{\ar@{->} "32"; "33"};
\endxy};
(20,0)*{T_{2,1}=0 \quad T_{2,3}=1};
(20,-5)*{T_{2,4}=0 \quad T_{2,5}=1};
(20,-10)*{T_{2,6}=0}
\endxy
\]

Next we define a function $S: \{1,\dots,n\} \lra \{+1,-1\}$ as follows.  First,

\begin{itemize}
\item For $1\leq i \leq n$, let $\sigma_i$ be the sign of crossing $i$; and
\item For $1\leq i \leq n$, let $d_i=+1$ if the ant goes \emph{over} crossing $i$ before it goes under crossing $i$, as it walks along the knot, and let $d_i=-1$ otherwise.
\end{itemize}

Now for $i=1,\dots,n$, let $S(i):= \sigma_id_i$.  We bundle this data into diagonal matrices, letting $\Sigma(K)$, $D(K)$ and $S(K)$ be the $\nxn$ diagonal matrices with the $\sigma_i$, the $d_i$ and the $S(i)$, respectively, along the diagonal,

\begin{equation*}
\Sigma(K):=Diag(\sigma_i), \quad D(K):=Diag(d_i), \quad S(K):= Diag(S(i))
\end{equation*}

Finally, for any indeterminate $x$, let $X^{\pm S}$ be the $\nxn$ diagonal matrix with entries $x^{\pm S(i)}$.

Then define

\begin{equation*}
\beta(K):= Det\big(I + T(I-X^{-S}) \big)
\end{equation*}

We illustrate these concepts with the sample long knot given above; we get the following:

\begin{equation*}
\Sigma(K)(x)=
\begin{pmatrix}
1 & 0 & 0 & 0 & 0 & 0 \\
0 & 1 & 0 & 0 & 0 & 0 \\
0 & 0 & 1 & 0 & 0 & 0 \\
0 & 0 & 0 & -1 & 0 & 0 \\
0 & 0 & 0 & 0 & -1 & 0 \\
0 & 0 & 0 & 0 & 0 & 1
\end{pmatrix}
\ \
D(K)(x)=
\begin{pmatrix}
-1 & 0 & 0 & 0 & 0 & 0 \\
0 & 1 & 0 & 0 & 0 & 0 \\
0 & 0 & -1 & 0 & 0 & 0 \\
0 & 0 & 0 & 1 & 0 & 0 \\
0 & 0 & 0 & 0 & -1 & 0 \\
0 & 0 & 0 & 0 & 0 & 1
\end{pmatrix}
\end{equation*}

\begin{equation*}
S(K)(x)=
\begin{pmatrix}
-1 & 0 & 0 & 0 & 0 & 0 \\
0 & 1 & 0 & 0 & 0 & 0 \\
0 & 0 & -1 & 0 & 0 & 0 \\
0 & 0 & 0 & -1 & 0 & 0 \\
0 & 0 & 0 & 0 & 1 & 0 \\
0 & 0 & 0 & 0 & 0 & 1
\end{pmatrix}
\ \
X^{-S}=
\begin{pmatrix}
x & 0 & 0 & 0 & 0 & 0 \\
0 & 1/x & 0 & 0 & 0 & 0 \\
0 & 0 & x & 0 & 0 & 0 \\
0 & 0 & 0 & x & 0 & 0 \\
0 & 0 & 0 & 0 & 1/x & 0 \\
0 & 0 & 0 & 0 & 0 & 1/x
\end{pmatrix}
\end{equation*}

\begin{equation*}
T(K)(x)=
\begin{pmatrix}
0 & 0 & 1 & 0 & 1 & 0 \\
0 & 0 & 1 & 0 & 1 & 0 \\
0 & 1 & 0 & 0 & 1 & 0 \\
0 & 1 & 0 & 0 & 1 & 1 \\
0 & 1 & 0 & 0 & 0 & 1 \\
1 & 1 & 1 & 0 & 1 & 0
\end{pmatrix}
\end{equation*}

\begin{equation*}
1+T(1-X^{-S})=
\begin{pmatrix}
1 & 0 & 1-x & 0 & 1-1/x & 0 \\
0 & 1 & 1-x & 0 & 1-1/x & 0 \\
0 & 1-1/x & 1 & 0 & 1-1/x & 0 \\
0 & 1-1/x & 0 & 1 & 1-1/x & 1-1/x \\
0 & 1-1/x & 0 & 0 & 1 & 1-1/x \\
1-x & 1-1/x & 1-x & 0 & 1-1/x & 1
\end{pmatrix}
\end{equation*}

One obtains $\beta(K)(x)= -\frac{1}{x^2}+ \frac{3}{x}-3 + 3x - x^2$.

The purpose of this paper is to prove the following theorem:

\begin{theorem}
\begin{equation*}
\beta(K)(x)=\pm x^{-l} \Delta(K)(x)
\end{equation*}
where $l$ is the number of crossings $j$ for which $S(j)=+1$.
\end{theorem}

\subsection{The Alexander Invariant}

We now recall the Alxander procedure, as adapted for long knots.  Once again, we assume the crossings and regions are numbered as per our convention (see Section \ref{Conventions}).

\subsubsection{Alexander's Procedure}

The first step in determining the Alexander invariant is to mark the crossings in any diagram of the knot according to the following scheme (for later convenience in our proof we have varied the usual marking scheme, but as explained below the result is equivalent to using the usual scheme):

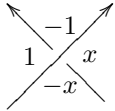
\[
\xy
(-20,0)*{
\xy
(-7,7)*{} = "1";
(7,7)*{} = "2";
(-7,-7)*{} = "3";
(7,-7)*+{} = "4";
(1,-1)*{} = "5";
(-1,1)*{} = "6";
(0,-4)*{-x};
(-4,0)*{1};
(0,4)*{-1};
(4,0)*{x};
{\ar@{-} "4"; "5"};
{\ar@{->} "6"; "1"};
{\ar@{->} "3"; "2"};
\endxy};
(20,0)*{}
\endxy
\]

In other words, as we proceed along the upper strand in the direction of the orientation, the $x$'s are on the right (first $-x$, then $x$), and the $1$'s are on the left (first $1$, then $-1$).

We now create an $\nxn$ incidence matrix, $A(K)$, whose $(i,j)$ entry contains the marking for crossing $i$ which lies in region $j$ (and $0$ if crossing $i$ and region $j$ do not abut).  The Alexander polynomial $\Delta(K)(x)$ is the determinant of this matrix, and is an invariant of the knot up to a factor of the form $\pm x^k$ for $k \in \Z$.

\subsubsection{An Example}

We illustrate the above marking procedure with the same knot which was used to illustrate the procedure of Bar-Natan and Dancso.  Vertex numbers are in parentheses, while region numbers are in square brackets.

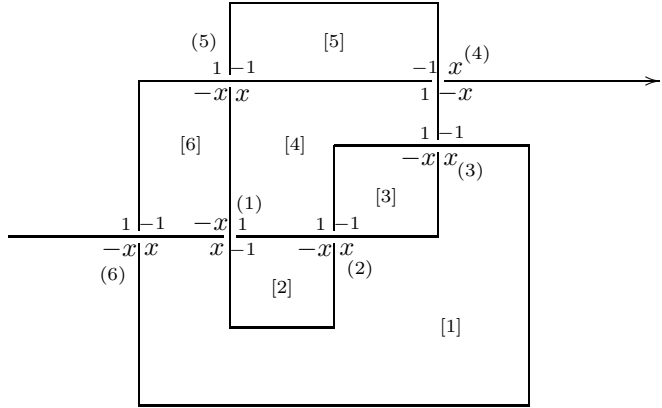
\[
\xy 0;/r.41pc/:
(-30,0)*{} = "1";
(-13.5,0)*{} = "2";
(-12.5,0)*{} = "3";
(3,0)*{} = "4";
(-13,-7)*{} = "5";
(-13,18)*{} = "6";
(3,18)*{} = "7";
(3,7.5)*{} = "8";
(3,6.5)*{} = "9";
(-5,-7)*{} = "11";
(-5,-0.5)*{} = "12";
(-5,0.5)*{} = "13";
(-5,7)*{} = "14";
(10,7)*{} = "15";
(-4,1)*{\scriptstyle{}} = "16";
(1,1)*{} = "17";
(1,16)*{} = "18";
(-14,16)*{} = "19";
(-14,-4)*{} = "20";
(-6,-4)*{} = "21";
(-6,-1)*{\scriptstyle{}} = "22";
(4,-7)*{\scriptstyle{[1]}};
(-1,3)*{\scriptstyle{[3]}};
(-8,7)*{\scriptstyle{[4]}};
(-5,15)*{\scriptstyle{[5]}};
(-16,7)*{\scriptstyle{[6]}};
(-9,-4)*{\scriptstyle{[2]}};
(-11.5,2.5)*{\scriptstyle{(1)}};
(-12,1)*{\scriptstyle{1}};
(-12,-1)*{\scriptstyle{-1}};
(-14.5,1)*{-x};
(-14,-1)*{x};
(-3,-2.5)*{\scriptstyle{(2)}};
(-4,-1)*{x};
(-4,1)*{\scriptstyle{-1}};
(-6.5,-1)*{-x};
(-6,1)*{\scriptstyle{1}};
(5.5,5)*{\scriptstyle{(3)}};
(4,8)*{\scriptstyle{-1}};
(2,8)*{\scriptstyle{1}};
(4,6)*{x};
(1.5,6)*{-x};
(6,14)*{\scriptstyle{(4)}};
(4,13)*{\ x};
(2,13)*{\scriptstyle{-1}};
(4,11)*{\ -x};
(2,11)*{\scriptstyle{1}};
(-15,15)*{\scriptstyle{(5)}};
(-12,13)*{\scriptstyle{-1}};
(-12,10.9)*{x};
(-14,13)*{\scriptstyle{1}};
(-14.5,11)*{-x};
(-22,-3)*{\scriptstyle{(6)}};
(-19,1)*{\scriptstyle{-1}};
(-21,1)*{\scriptstyle{1}};
(-21.5,-1)*{-x};
(-19,-1)*{x};
(10,5)*{} = "23";
(10,-13)*{} = "24";
(-20,-13)*{} = "25";
(-20,-0.5)*{} = "26";
(-20,0.5)*{} = "27";
(-20,12)*{} = "28";
(-13,11.5)*{} = "29";
(-13,12.5)*{} = "30";
(2.5,12)*{} = "31";
(3.5,12)*{} = "32";
(20,12)*{} = "33";
{\ar@{-} "1"; "2"};
{\ar@{-} "3"; "4"};
{\ar@{-} "4"; "9"};
{\ar@{-} "8"; "7"};
{\ar@{-} "7"; "6"};
{\ar@{-} "6"; "30"};
{\ar@{-} "5"; "11"};
{\ar@{-} "11"; "12"};
{\ar@{-} "13"; "14"};
{\ar@{-} "14"; "15"};
{\ar@{-} "15"; "24"};
{\ar@{-} "24"; "25"};
{\ar@{-} "25"; "26"};
{\ar@{-} "27"; "28"};
{\ar@{-} "28"; "31"};
{\ar@{-} "29"; "5"};
{\ar@{->} "32"; "33"};
\endxy
\]

This leads to the following Alexander matrix:

\begin{equation*}
A(K)(x)=
\begin{pmatrix}
x & -1 & 0 & 1 & 0 & -x \\
x & -x & -1 & 1 & 0 & 0 \\
x & 0 & -x & 1 & 0 & 0 \\
0 & 0 & 0 & 1 & -1 & 0 \\
0 & 0 & 0 & x & -1 & -x \\
x & 0 & 0 & 0 & 0 & -1
\end{pmatrix}
\end{equation*}

The determinant of this matrix is $\Delta(K)(x)= x-3x^2+3x^3-3x^4+x^5$.  It is clear that, in this particular example at least, $\beta(K)(x)= -x^{-3} \Delta(K)(x)$, and 3 is the number of crossings for which $S(j)=+1$.

\subsubsection{Remarks on our Conventions}

As indicated earlier, Alexander's paper \cite{Alex} dealt not with long knots but round knots; we have adapted the construction accordingly.  Also, Alexander's procedure involved numbering all $(n+2)$ regions and preparing an $n \times (n+2)$ matrix as above, and then deleting the two columns corresponding to any two adjacent regions.  Alexander showed that the resulting polynomial was independent of the choice of columns to delete.  The construction we have given corresponds to deleting the columns relating to the two unbounded regions.

Again as indicated earlier, we have also changed the usual marking convention (the one adopted in \cite{Alex}), which follows the scheme:

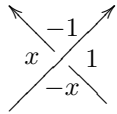
\[
\xy
(-20,0)*{
\xy
(-7,7)*{} = "1";
(7,7)*{} = "2";
(-7,-7)*{} = "3";
(7,-7)*+{} = "4";
(1,-1)*{} = "5";
(-1,1)*{} = "6";
(0,-4)*{-x};
(-4,0)*{x};
(0,4)*{-1};
(4,0)*{1};
{\ar@{-} "4"; "5"};
{\ar@{->} "6"; "1"};
{\ar@{->} "3"; "2"};
\endxy};
(20,0)*{}
\endxy
\]

In other words, as we proceed along the \emph{lower} strand in the direction of the orientation, the $x$'s are on the left (first $-x$, then $x$), and the $1$'s are on the right (first $1$, then $-1$).  However, the Alexander polynomial is invariant under rotation of the long knot by $\pi$ radians along its axis.  Such a rotation flips top strands into bottom strands (and right into left) and vice versa, and therefore carries one marking scheme into the other.  By invariance, the specific choice of convention does not change the result.

\section{Proof of the Conjecture}

Again we assume given a long knot $K$ with numbered crossings and regions, as per our convention (see Section \ref{Conventions}).

We will replace the Alexander matrix $A(K)$ of the knot by a new matrix obtained as follows.

For any region $i$ and crossing $j$ of $K$ we define the winding number $w(i,j)$ to be the winding number of any point in region $i$ with respect to $Sp(j)$.  Now define the $\nxn$ `winding' matrix $W(K)$ of $K$ to be:

\begin{equation}
\label{WindingMatrix}
W(K) := \big( w(i,j) \big)
\end{equation}

Note that to compute the winding number of a region with respect to $Sp(j)$, we can pick any point within the region, draw any ray from that point to infinity, and count the number of times the ray crosses $Sp(j)$.  We count $+1$ if $Sp(j)$ crosses from left to right as we move from the point out to infinity, and $-1$ if $Sp(j)$ crosses from right to left.  We can arrange that all crossings between $Sp(j)$ and the ray are transversal by deforming $Sp(j)$ slightly, if needed.

We illustrate the computation of this winding matrix with the same sample long knot as previously.  The following shows the span $Sp(2)$, and rays which may be used to determine the winding number of each region with respect to that particular span.  As before, vertex numbers are in parentheses and region number are in square brackets.  We have indicated the winding number of each region (with respect to $Sp(2)$) by a number in curly brackets (e.g. $\{-1\}$) at the end of the ray used to compute that number.

\[
\xy 0;/r.41pc/:
(-30,0)*{} = "1";
(-13.5,0)*{} = "2";
(-12.5,0)*{} = "3";
(3,0)*{} = "4";
(-13,-7)*{} = "5";
(-13,18)*{} = "6";
(3,18)*{} = "7";
(3,7.5)*{} = "8";
(3,6.5)*{} = "9";
(-5,-7)*{} = "11";
(-5,-0.5)*{} = "12";
(-5,0.5)*{} = "13";
(-5,7)*{} = "14";
(10,7)*{} = "15";
(-4,1)*{\scriptstyle{}} = "16";
(1,1)*{} = "17";
(1,16)*{} = "18";
(-14,16)*{} = "19";
(-14,-4)*{} = "20";
(-6,-4)*{} = "21";
(-6,-1)*{\scriptstyle{}} = "22";
(4,-7)*{\scriptstyle{[1]}};
(-2,3)*{\scriptstyle{[3]}};
(-8,7)*{\scriptstyle{[4]}};
(-5,15)*{\scriptstyle{[5]}};
(-16,7)*{\scriptstyle{[6]}};
(-9,-2)*{\scriptstyle{[2]}};
(-11.5,1)*{\scriptstyle{(1)}};
(-3.5,-1.5)*{\scriptstyle{(2)}};
(-4,1)*{*};
(-6,-1.5)*{*};
(5,6)*{\scriptstyle{(3)}};
(5,13)*{\scriptstyle{(4)}};
(-14.5,13)*{\scriptstyle{(5)}};
(-21.5,-1.5)*{\scriptstyle{(6)}};
(10,5)*{} = "23";
(10,-13)*{} = "24";
(-20,-13)*{} = "25";
(-20,-0.5)*{} = "26";
(-20,0.5)*{} = "27";
(-20,12)*{} = "28";
(-13,11.5)*{} = "29";
(-13,12.5)*{} = "30";
(2.5,12)*{} = "31";
(3.5,12)*{} = "32";
(20,12)*{} = "33";
{\ar@{.} (-4,1); (2,1)};
{\ar@{.} (2,1); (2,19)};
{\ar@{.} (2,19); (-14,19)};
{\ar@{.} (-14,19); (-14,-8)};
{\ar@{.} (-14,-8); (-6,-8)};
{\ar@{.} (-6,-8); (-6,-1.5)};
{\ar@{-} "1"; "2"};
{\ar@{-} "3"; "4"};
{\ar@{-} "4"; "9"};
{\ar@{-} "8"; "7"};
{\ar@{-} "7"; "6"};
{\ar@{-} "6"; "30"};
{\ar@{-} "5"; "11"};
{\ar@{-} "11"; "12"};
{\ar@{-} "13"; "14"};
{\ar@{-} "14"; "15"};
{\ar@{-} "15"; "24"};
{\ar@{-} "24"; "25"};
{\ar@{-} "25"; "26"};
{\ar@{-} "27"; "28"};
{\ar@{-} "28"; "31"};
{\ar@{-} "29"; "5"};
{\ar@{->} "32"; "33"};
(0,3)*{\bullet};
{\ar@{~>} (0,3); (20,3)};
(22,3)*{\scriptstyle{\{-1\}}};
(0,9)*{\bullet};
{\ar@{~>} (0,9); (20,9)};
(22,9)*{\scriptstyle{\{-1\}}};
(-16,9)*{\bullet};
{\ar@{~>} (-16,9); (-30,9)};
(-32,9)*{\scriptstyle{\{0\}}};
(0,15)*{\bullet};
{\ar@{~>} (0,15); (20,15)};
(22,15)*{\scriptstyle{\{-1\}}};
(2,-3)*{\bullet};
{\ar@{~>} (2,-3); (20,-3)};
(22,-3)*{\scriptstyle{\{0\}}};
(-9,-4)*{\bullet};
{\ar@{~>} (-9,-4); (-30,-4)};
(-32,-4)*{\scriptstyle{\{-1\}}};
\endxy
\]

We get the following winding matrix:

\begin{equation*}
W(K)(x)=
\begin{pmatrix}
0 & 0 & 0 & 1 & 1 & 1 \\
0 & -1 & -1 & 0 & 0 & 0 \\
-1 & -1 & 0 & 1 & 1 & 0 \\
-1 & -1 & -1 & 0 & 0 & -1 \\
-1 & -1 & -1 & -1 & 0 & -1 \\
0 & 0 & 0 & 1 & 1 & 0
\end{pmatrix}
\end{equation*}

Lastly, we let $X^{\pm (1+S)/2}$ be the diagonal matrix whose $i$th diagonal is $x^{\pm 1}$ if $S(i)=+1$, and $x^0=1$ if $S(i)=-1$.

We will show that

\begin{proposition}
\begin{equation}
\label{B(K)Formula}
\Sigma(K) A(K) W(K) S(K) X^{-(1+S)/2} = \big( 1+T^t(1-X^{-S})\big)
\end{equation}
where the superscript in $T^t$ means we take the transpose of the matrix $T$.
\end{proposition}

In our ongoing example, we get the following:

\begin{equation*}
X^{-(1+S)/2}=
\begin{pmatrix}
1 & 0 & 0 & 0 & 0 & 0 \\
0 & 1/x & 0 & 0 & 0 & 0 \\
0 & 0 & 1 & 0 & 0 & 0 \\
0 & 0 & 0 & 1 & 0 & 0 \\
0 & 0 & 0 & 0 & 1/x & 0 \\
0 & 0 & 0 & 0 & 0 & 1/x
\end{pmatrix}
\end{equation*}

\begin{equation*} 
1+T^t(1-X^{-S})=
\begin{pmatrix}
1 & 0 & 0 & 0 & 0 & 1-1/x \\
0 & 1 & 1-x & 1-x & 1-1/x & 1-1/x \\
1-x & 1-1/x & 1 & 0 & 0 & 1-1/x \\
0 & 0 & 0 & 1 & 0 & 0 \\
1-x & 1-1/x & 1-x & 1-x & 1 & 1-1/x \\
0 & 0 & 0 & 1-x & 1-1/x & 1
\end{pmatrix}
\end{equation*}

One can verify that, indeed, Equation (\ref{B(K)Formula}) holds in this case.  Moreover, the determinant of
both sides is $-\frac{1}{x^2}+ \frac{3}{x}-3 + 3x - x^2$, that is to say $\beta(K)(x)$.

To prove the Proposition we simply compare the entries on both sides and find that they are the same.  We start with the diagonal entries.  We will first prove the following:

\begin{lemma}
$\big(A(K)W(K)\big)_{jj}= \big(D(K) X^{(1+S)/2} \big)_{jj}$.
\end{lemma}

Entry $(j,j)$ of $A(K)W(K)$ is $\sum_l A(K)_{jl} W(K)_{lj}$.  Since $A(K)_{jl}$ is just the marking of crossing $j$ in region $l$, there are (up to) four non-zero summands, corresponding to the four regions abutting crossing $j$.  In general, some of these regions may actually be non-numbered (i.e., may be unbounded regions), and hence will not contribute to this sum, since $A(K)$ only includes markings in numbered (i.e., bounded) regions.  However, for purposes of proving the above Lemma (and also the next Lemma) it will be useful to think of the sum  $\sum_l A(K)_{jl} W(K)_{lj}$ as including summands for all four regions abutting a crossing.  Since it is clear that the unbounded regions have winding number zero with respect to the span of any crossing, there is no harm in thinking of the above sums as also including summands $A(K)_{j,u} w(u,j)$ (where $A(K)_{j,u}$ denotes the marking of crossing $j$ in an unbounded region $u$, and $w(u,j)$ is the winding number of region $u$ with respect to $Sp(j)$), for any unbounded region which abuts the crossing.

It will be useful to draw pictures of the crossings in order to clarify the argument (a sample picture appears below).  If we consider the four regions abutting crossing $j$, only one of these touches $Sp(j)$ both before and after crossing $j$; we call this region the Spanning Region, and we call the remaining three regions the Complementary Regions.  In our diagrams, a bullet at the tail of a strand leading into crossing $j$ indicates along which strand an ant traveling along the knot first approaches the crossing.  The portion of $Sp(j)$ outside a small neighbourhood of $j$ is depicted with loosely spaced dots, to indicate that we make no attempt to draw that part of the diagram precisely.  For instance, the following picture shows a crossing with $\sigma_j =d_j = 1$:

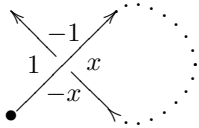
\[
\xy
(-7,7)*{} = "1";
(7,7)*{} = "2";
(-7,-7)*{\bullet} = "3";
(7,-7)*+{} = "4";
(1,-1)*{} = "5";
(-1,1)*{} = "6";
(0,-4)*{-x};
(-4,0)*{1};
(0,4)*{-1};
(4,0)*{x};
{\ar@{>-} "4"; "5"};
{\ar@{->} "6"; "1"};
{\ar@{->} "3"; "2"};
"2"; "4"**\crv{~*=<2mm>{.}(10,10) & (25,0) & (10,-10)};
\endxy
\]

\begin{proof}

We consider the four cases $\sigma_j=\pm1$ and $d_j=\pm1$ in turn.

We start with $\sigma_j=d_j=S(j)=1$.  We claim that the winding number is exactly one greater in the Spanning Region than in the Complementary Regions.  To compute the winding numbers in the Spanning Region and any of the Complementary Regions, we pick a point A in the Spanning Region and a point B in the Complementary Region, close enough together that the straight line segment $(A,B)$ crosses $Sp(j)$ exactly once, and does not otherwise cross the knot.  We draw a ray from B through A and in a straight line on to infinity.  We compute $\w(B,j)$ using this ray, and we compute $\w(A,j)$ using the portion of the ray complementary to the line segment $(A,B)$:

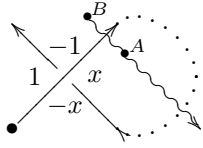
\[
\xy
(-7,7)*{} = "1";
(7,7)*{} = "2";
(-7,-7)*{\bullet} = "3";
(7,-7)*{} = "4";
(1,-1)*{} = "5";
(-1,1)*{} = "6";
(0,-4)*{-x};
(-4,0)*{1};
(0,4)*{-1};
(4,0)*{x};
(9,4)*{\scriptstyle{\ A}};
(4,9)*{\scriptstyle{\ B}};
(8,3)*{\scriptstyle{\bullet}} = "10";
(3,8)*{\scriptstyle{\bullet}} = "11";
{\ar@{>-} "4"; "5"};
{\ar@{->} "6"; "1"};
{\ar@{->} "3"; "2"};
{\ar@{~>} (3,8); (18,-7)};
"2"; "4"**\crv{~*=<2mm>{.}(10,10) & (25,0) & (10,-10)};
\endxy
\]

It is now clear that $w(B,j)=w(A,j)-1$.  Hence we get the following for $\big(A(K)W(K)\big)_{jj}$:

\begin{equation*}
(-1+1-x)(w(A,j)-1) + x. w(A,j) = x = d_j x^{(1+S(j))/2}
\end{equation*}
as claimed.

The remaining cases are summarized in the following diagrams:

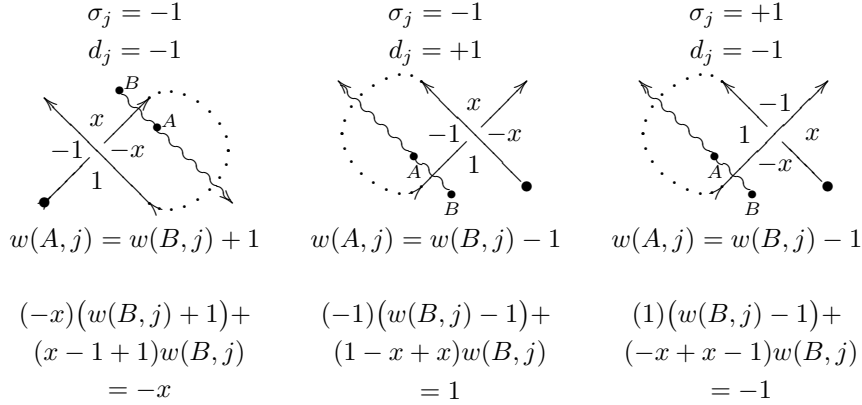
\[
\xy
(-40,15)*{\sigma_j=-1};
(-40,10)*{d_j=-1};
(-40,-2)*{
\xy
(-7,7)*{} = "1";
(7,7)*{} = "2";
(-7,-7)*{\bullet} = "3";
(7,-7)*{} = "4";
(-1,-1)*{} = "5";
(1,1)*{} = "6";
(0,-4)*{1};
(-4,0)*{-1};
(0,4)*{x};
(4,0)*{-x};
{\ar@{>->} "4"; "1"};
{\ar@{->} "6"; "2"};
{\ar@{>-} "3"; "5"};
(9,4)*{\scriptstyle{\ A}};
(4,9)*{\scriptstyle{\ B}};
(8,3)*{\scriptstyle{\bullet}} = "10";
(3,8)*{\scriptstyle{\bullet}} = "11";
{\ar@{~>} (3,8); (18,-7)};
"2"; "4"**\crv{~*=<2mm>{.}(10,10) & (25,0) & (10,-10)};
\endxy};
(-40,-15)*{\w(A,j)=\w(B,j)+1};
(-40,-25)*{(-x)\big(\w(B,j)+1\big) +};
(-40,-30)*{\ (x-1+1)\w(B,j)};
(-40,-35)*{\ =-x};
(0,15)*{\sigma_j=-1};
(0,10)*{d_j=+1};
(0,-2)*{
\xy
(-7,7)*{} = "1";
(7,7)*{} = "2";
(-7,-7)*{} = "3";
(7,-7)*{\bullet} = "4";
(-1,-1)*{} = "5";
(1,1)*{} = "6";
(0,-4)*{1};
(-4,0)*{-1};
(0,4)*{x};
(4,0)*{-x};
{\ar@{->} "4"; "1"};
{\ar@{->} "6"; "2"};
{\ar@{>-} "3"; "5"};
"1"; "3"**\crv{~*=<2mm>{.}(-10,10) & (-25,0) & (-10,-10)};
(-3,-8)*{\scriptstyle{\bullet}} = "7";
(-8,-3)*{\scriptstyle{\bullet}} = "8";
(-8,-5)*{\scriptstyle{A}};
(-3,-10)*{\scriptstyle{B}};
{\ar@{~>} (-3,-8); (-18,7)};
\endxy};
(0,-15)*{\w(A,j)=\w(B,j)-1};
(0,-25)*{(-1)\big(\w(B,j)-1\big) +};
(0,-30)*{\ (1-x+x)\w(B,j)};
(0,-35)*{\ =1};
(40,15)*{\sigma_j=+1};
(40,10)*{d_j=-1};
(40,-2)*{
\xy
(-7,7)*{} = "1";
(7,7)*{} = "2";
(-7,-7)*{} = "3";
(7,-7)*{\bullet} = "4";
(-1,1)*{} = "5";
(1,-1)*{} = "6";
(0,-4)*{-x};
(-4,0)*{1};
(0,4)*{-1};
(5,0)*{x};
(-3,-8)*{\scriptstyle{\bullet}} = "7";
(-8,-3)*{\scriptstyle{\bullet}} = "8";
(-8,-5)*{\scriptstyle{A}};
(-3,-10)*{\scriptstyle{B}};
{\ar@{->} "5"; "1"};
{\ar@{-} "4"; "6"};
{\ar@{>->} "3"; "2"};
"1"; "3"**\crv{~*=<2mm>{.}(-10,10) & (-25,0) & (-10,-10)};
{\ar@{~>} (-3,-8); (-18,7)};
\endxy};
(40,-15)*{\w(A,j)=\w(B,j)-1};
(40,-25)*{(1)\big(\w(B,j)-1\big) +};
(40,-30)*{\ (-x+x-1)\w(B,j)};
(40,-35)*{\ =-1}
\endxy
\]
In each case it is clear that we obtained $d_jx^{(1+S(j))/2}$, as claimed.
\end{proof}

We now consider the off-diagonal entries of $A(K)W(K)$.  We will prove that, for $1\leq i \ne j \leq n$:

\begin{lemma}
$\big(A(K)W(K)\big)_{ij} = \sigma_i (x-1) T_{ji}$.
\end{lemma}

Entry $\big(A(K)W(K)\big)_{ij}$ is the sum $\sum_l A(K)_{il}W(K)_{lj}$, which may have up to four non-zero summands, one in respect of each numbered region abutting crossing $i$.  As with the previous Lemma, it will be useful to think of this sum as including summands for all four regions abutting a crossing, including any of these which are non-numbered, if necessary by `adding' summands $A(K)_{j,u} w(u,j)$ in respect of any non-numbered regions which abut crossing $i$ (such summands being in any case equal to zero).

\begin{proof}
If $T_{ji}=0$, then either $i$ is not in $Sp(j)$, or $Sp(j)$ goes over, but not under, crossing $i$.  In the former case, all four regions around crossing $i$ have the same winding number relative to $Sp(j)$, so their coefficients with respect to crossing $i$ cancel in computing $\big(A(K)W(K)\big)_{ij}$.  In the latter case, we have the following picture (if crossing $i$ is negative):

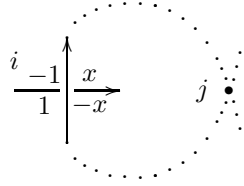
\[
\xy
(-12,0)*{} = "1";
(-6,0)*{} = "2";
(-4,0)*{} = "3";
(2,0)*{} = "4";
(-5,7)*{} = "5";
(-5,-7)*{} = "6";
(-8,2)*{-1};
(-8,-2)*{1};
(-2,-2)*{-x};
(-2,2)*{x};
(-12,4)*{i};
(13,0)*{j} ="7";
(16.5,0)*{\scriptstyle{\bullet}} ="7";
{\ar@{-} "1"; "2"};
{\ar@{->} "3"; "4"};
{\ar@{->} "6"; "5"};
"5"; (18,-5)**\crv{~*=<2mm>{.}(-5,10) & (8,15) & (16,5)};
"6"; (18,5)**\crv{~*=<2mm>{.}(-5,-10) & (8,-15) & (16,-5)};
\endxy
\]
where $Sp(j)$ follows the vertical strand over crossing $i$ and then continues along the dotted line.  The regions on the same side of the vertical strand have the same winding number, and the resulting contributions have the factor $(x-x)=0$ or $(1-1)=0$, so there is no contribution (the story is essentially the same for a positive crossing).

We now suppose $T_{ji}=1$.  Then the winding numbers for the two regions to the left of the under strand (as we go through crossing $i$ in the direction of the orientation) are equal, and are exactly one less than the winding numbers for the regions to the right.  This can be seen from the following diagram, which shows the computation $\big(A(K)W(K)\big)_{ij} = \sigma_i (x-1)$ when $i$ is a positive crossing (left picture) or a negative crossing (right picture).

\[
\xy
(-25,0)*{
\xy
(-12,0)*{} = "1";
(-5,-1)*{} = "2";
(-5,1)*{} = "3";
(2,0)*{} = "4";
(-5,7)*{} = "5";
(-5,-7)*{} = "6";
(-8,2)*{1};
(-8,-2)*{-x};
(-2,-2)*{x};
(-2,2)*{-1};
(-13,5)*{i};
(13,0)*{j};
(16.5,0)*{\scriptstyle{\bullet}};
{\ar@{->} "1"; "4"};
{\ar@{-} "6"; "2"};
{\ar@{->} "3"; "5"};
"5"; (18,-5)**\crv{~*=<2mm>{.}(-5,10) & (8,15) & (16,5)};
"6"; (18,5)**\crv{~*=<2mm>{.}(-5,-10) & (8,-15) & (16,-5)};
(-8,5)*{\scriptstyle{\bullet}} = "7";
(4,5)*{\scriptstyle{\bullet}} = "8";
(-8,7)*{\scriptstyle{A}};
(4,7)*{\scriptstyle{B}};
{\ar@{~>} (-8,5);(23,5)};
\endxy};
(-25,-15)*{\w(A,j)=\w(B,j)-1};
(-25,-25)*{(1-x)\big(\w(B,j)-1\big) +};
(-25,-30)*{\ (x-1)\w(B,j)};
(-25,-35)*{\ =x-1};
(25,0)*{
\xy
(-12,0)*{} = "1";
(-5,-1)*{} = "2";
(-5,1)*{} = "3";
(2,0)*{} = "4";
(-5,7)*{} = "5";
(-5,-7)*{} = "6";
(-8,2)*{x};
(-8,-2)*{-1};
(-2,-2)*{1};
(-2,2)*{-x};
(-13,5)*{i};
(13,0)*{j};
(16.5,0)*{\scriptstyle{\bullet}};
{\ar@{->} "4"; "1"};
{\ar@{->} "3"; "5"};
{\ar@{-} "6"; "2"};
"5"; (18,-5)**\crv{~*=<2mm>{.}(-5,10) & (8,15) & (16,5)};
"6"; (18,5)**\crv{~*=<2mm>{.}(-5,-10) & (8,-15) & (16,-5)};
(-8,5)*{\scriptstyle{\bullet}} = "7";
(4,5)*{\scriptstyle{\bullet}} = "8";
(-8,7)*{\scriptstyle{A}};
(4,7)*{\scriptstyle{B}};
{\ar@{~>} (-8,5);(23,5)};
\endxy};
(25,-15)*{\w(A,j)=\w(B,j)-1};
(25,-25)*{(x-1)\big(\w(B,j)-1\big) +};
(25,-30)*{\ (1-x)\w(B,j)};
(25,-35)*{\ =1-x}
\endxy
\]
\end{proof}

\begin{proof}[Proof of Proposition]

We now add in the effect of left multiplying $A(K)W(K)$ by $\Sigma(K)$ and right multiplying by $S(K)X^{-(1+S)/2}$.  For the $i$th diagonal entry, we get:
\begin{equation*}
\sigma_i \big(A(K)W(K)\big)_{ii} S(i) x^{-(1+S(i))/2} = \sigma_i (d_i x^{(1+S(i))/2}) (\sigma_id_i) x^{-(1+S(i))/2} = 1
\end{equation*}

For the $(i,j)$th off-diagonal entry, we get:
\begin{align*}
\sigma_i \big(A(K)W(K)\big)_{ij} S(j) x^{-(1+S(j))/2} & = \sigma_i (\sigma_i (x-1)T_{ji}) S(j) x^{-(1+S(j))/2} \\
& = (x-1) T_{ji} S(j) x^{-(1+S(j))/2}
\end{align*}

If $S(j)=+1$, we get
\begin{equation*}
(x-1) T_{ji} S(j) x^{-(1+S(j))/2} = (x-1)T_{ji} .1.x^{-1}=(1-1/x)T_{ji} =(1-x^{-S(j)})T_{ji}
\end{equation*}

Similarly, if $S(j)=-1$, we get
\begin{equation*}
(x-1) T_{ji} S(j) x^{-(1+S(j))/2} = (x-1)T_{ji} (-1)(1)=(1-x)T_{ji} =(1-x^{-S(j)})T_{ji}
\end{equation*}

Putting all this together, we find
\begin{equation*}
\Sigma(K) A(K) W(K) S(K) X^{-(1+S)/2}=\big( 1 + T^t(1-X^{-S})\big)
\end{equation*}
as required.

\end{proof}

\begin{proof}[Proof of the Theorem]

\begin{align*}
\beta(K)(x)= Det(1+T(1-X^{-S})) &= Det((1-X^{-S})^{-1} +T) Det(1-X^{-S}) \\
&=  Det(1-X^{-S}) Det((1-X^{-S})^{-1} +T) \\
&= Det(1+(1-X^{-S})T) \\
&= Det((1+(1-X^{-S})T)^t) \\
&= Det(1+T^t (1-X^{-S})) \\
&= Det(\Sigma(K) A(K) W(K) S(K) X^{-(1+S)/2}) \\
&= \pm \Delta(K)(x) Det(W(K)) x^{-l}
\end{align*}
where $l$ is the number of crossings $j$ for which $S(j)=1$.

Since it was shown in \cite{DBN-WKO} that $Det(1+T(1-X^{-S}))$ is not zero, it follows that $Det(W(K)) \ne 0$. But since $W(K)$ has integer coefficients, it follows that $Det(W(K))=\pm 1$, and the theorem is proved.

\end{proof}

\end{document}